\documentclass{amsart} \textwidth=14.5cm \oddsidemargin=1cm
\usepackage{amsmath}
\usepackage{amsxtra}
\usepackage{amscd}
\usepackage{amsthm}
\usepackage{amsfonts}
\usepackage{amssymb}
\usepackage{eucal}

\newtheorem{Thm}{Theorem}
\newtheorem{Cor}{Corollary}
\newtheorem{Lem}{Lemma}
\newtheorem{Prop}{Proposition}

\theoremstyle{remark}
\newtheorem{Rem}{Remark}

\newcommand{\Ql}{{\overline {{\Bbb Q}_l} }}

\newcommand{\cal}{\mathcal}
\newcommand{\Gr}{{{\cal G}{\frak r} }}
\newcommand{\Fl}{{{\cal F}\ell}}

\newcommand{\bu}{\bullet}

\newcommand{\To}{\longrightarrow}
\newcommand{\iso}{{\widetilde \longrightarrow}}

\newcommand{\imbed}{\hookrightarrow}
\newcommand{\codim}{{\rm codim}\ }
\newcommand{\rank}{{\rm rank}\ }
\newcommand{\sur}{\twoheadrightarrow}
\newcommand{\<}{\langle}
\renewcommand{\>}{\rangle}

\def\square{\hbox{\vrule\vbox{\hrule\phantom{o}\hrule}\vrule}}
\newcommand{\epf}{\square}

\newcommand{\Cbar}{{\bar C}}

\newcommand{\fbar}{{\bar f}}

\renewcommand{\P}{{\cal P}}
\newcommand{\PI}{{\cal P}_I}
\newcommand{\fP}{{^f{\cal P}}}
\newcommand{\fPI}{{^f{\cal P}_I}}
\newcommand{\fPfI}{{^f{\cal P}^f_I}}

\newcommand{\fPf}{{^f{\cal P}^f}}

\newcommand{\DIW}{D_{{\cal {IW}}}}

\newcommand{\fPhi}{{^f{\Phi}}}

\newcommand{\fPhif}{{^f{\Phi}^f}}

\newcommand{\Xtil}{\tilde{X}}

\newcommand{\N}{{\cal N}}
\newcommand{\Ntil}{{\tilde{\cal N}}}

\newcommand{\Nt}{{\tilde{\cal N}}}

\newcommand{\St}{{St}}

\newcommand{\oplusl}{\bigoplus\limits}
\newcommand{\cupl}{\bigcup\limits}
\newcommand{\fW}{{^f W}}

\newcommand{\fWf}{{^fW^f}}

\newcommand{\M}{{\cal M}}

\newcommand{\gp}{{\frak p}}

\newcommand{\Lg}{{ \frak g\check{\ }}}

\newcommand{\LG}{{ G\check{\ }}}

\renewcommand{\L}{{\cal L}}

\renewcommand{\O}{{\cal O}}
\newcommand{\Ohat}{{\hat\O}}
\newcommand{\Obar}{{\overline O}}

\newcommand{\F}{{\cal F}}

\newcommand{\A}{{\cal A}}
\newcommand{\B}{{\cal B}}
\newcommand{\C}{{\cal C}}

\newcommand{\J}{{\cal J}}

\newcommand{\Z}{{\cal Z}}
\newcommand{\cJ}{{\cal J}}
\newcommand{\Zet}{{\Bbb Z}}

\newcommand{\bbone}{{\bf 1}}
\renewcommand{\c}{{\underline {c}}}
\newcommand{\LRleq}{\underset{LR}{\leq}}

\newcommand{\unO}{{\underline{\cal O}}}

\newcommand{\bI}{{\bf I}}

\newcommand{\GO}{{\bf G_O}}

\newcommand{\Db}{D^b}

\newcommand{\bbK}{{\mathbb K}}

\newcommand{\bbV}{{\mathbb V}}

\newcommand{\Pone}{{\mathbb P}^1}
\newcommand{\Pn}{{\mathbb P}^n}

\newcommand{\Sa}{S}
\newcommand{\orb}{{\bf O}}

\author{Roman Bezrukavnikov}
\title[Affine flags and  nilpotent cone]
{Perverse sheaves on affine flags and nilpotent cone of
the Langlands dual group}
\begin{document}
\maketitle
\tableofcontents

\begin{abstract} This paper is a continuation of
\cite{1}. In \cite{1} we constructed an equivalence between
the derived category of equivariant coherent sheaves on the
cotangent bundle to the flag variety of a simple algebraic group and
a (quotient of) a category of constructible sheaves on the
affine flag variety of the Langlands dual group. Below we
prove certain properties of this equivalence related to cells
in the affine Weyl group; provide a similar
``Langlands dual'' description for the category of equivariant
coherent sheaves on the nilpotent cone, and link it to 
perverse coherent sheaves; and deduce some conjectures
by Lusztig and Ostrik.
%As an application we prove
%some conjectures about the asymptotic Hecke algebra, and get some
%new properties of the central sheaves introduced in \cite{KGB}.
\end{abstract}

\thanks{{\bf Acknowledgements.}
I am greatful to all the people mentioned in
 acknowledgements in \cite{1}. 
I also thank Eric Sommers for stimulating interest and the referee for
 helpful comments.
 The author worked on this paper while
supported by NSF  grant DMS-0071967  and employed by the 
Clay Institute.}

\begin{section}{Statements}
\subsection{Recollection of notations and set-up}
We keep the set-up and  notations of \cite{1}. In particular,
$\Fl$ is the affine flag variety  of a split simple group $G$
over an algebraically closed field $k$;           %=\Fpbar$ 
%which is either finite or algebraically closed; 
$W_f$ is the Weyl group of $G$,
and $W$ is the extended affine Weyl group;
$\fWf\subset \fW\subset W$ are the sets of minimal length
representatives of respectively 2-sided and left cosets
of $W_f$ in  $W$; $D_I=D_I(\Fl)$ is the Iwahori equivariant derived category of $l$-adic sheaves
($l\ne char (k)$) on $\Fl$ and
$\PI\subset D_I(\Fl)$ is the full subcategory of 
perverse sheaves. 
%(the notation used only for an algebraically
%$closed $k$); 
%while $\PI^{mix}$ is the category whose objects
%are mixed Iwahori equivariant perverse sheaves on $\Fl$, and morphisms
%are weight 0 geometric morphisms, i.e. weight 0
% morphisms between the pull-backs of sheaves to $\Fl_{\kbar}$
%(the notation is used for finite $k$ only).

$L_w$, $w\in W$ are irreducible objects of $\PI$.
%, or irreducible self-dual
%objects of $\PI^{mix}$ depending on the cardinality of $k$.
The Serre quotient category $\fPI$ %, $\fPI^{mix}$ 
of $\PI$ %, $\PI^{mix}$ are
is defined by
 $$\fP_I=\P_I/\< L_w\ |\ w\not \in \fW\>,$$
%$$\fP_I^{mix}=\P_I^{mix}/\< L_w(i)\ |\ w\not \in \fW, \ i\in \Zet\>,$$
where for  an abelian category $\A$, and a set $S$ of irreducible
objects of $\A$ we let $\<S\rangle$
denote the full abelian subcategory of objects obtained from elements
of $S$ by extensions.

$\N$ is the variety of nilpotent elements in the Langlands dual
Lie algebra $\Lg$; and $p_{Spr}:\Nt\to \N$ is its Springer
resolution. For an algebraic group $H$ acting on a variety $X$ we
write $D^H(X)$ instead of $D^b(Coh^H(X))$.

Convolution provides $D_I(\Fl)$ with a monoidal structure. 
In \cite{1} we constructed a monoidal functor
$$F:D^\LG(\Ntil)\to D_I(\Fl);$$
%$$F^{mix}:D^{\LG\times \Gm}(\Ntil)\to D_I^{mix}(\Fl);$$
and used it to define  an equivalence
\begin{equation}\label{feq}
%\begin{array}{ll}&
\fPhi:D^\LG(\Ntil)\iso \Db(\fP_I).
%\\&\fPhi^{mix}:D^{\LG\times \Gm}(\Ntil)\iso \Db(\fP_I^{mix}).
%\end{array}
\end{equation}

We now state the results proved in this note.

\subsection{Category $\fPfI$ %, $\fPfI^{mix}$
 and the nil-cone}
Let us  define a further Serre quotient category $\fPfI$ %, $\fPfI^{mix}$
of $\PI$ by
 $$\fPf_I=\P_I/\< L_w\ |\ w\not \in \fWf\>.$$
%$$\fPfI^{mix}=\P_I^{mix}/\< L_w(i)\ |\ w\not \in \fWf, \ i\in \Zet\>.$$

Let $pr_{ff}:\fPI\to \fPfI$ %, $pr_{ff}^{mix}:\fPI^{mix}\to \fPfmI$
 be the projection functor. (We will use the same notation for the
extension of these exact functors to the derived categories.
We will also abuse notations
by omitting the projection to a quotient category functor from notations;
e.g. we will sometimes write ``$X$'' or ``$X$ considered as an object of
$\fPf$'' instead of $pr_{ff}(X)$.)

\begin{Thm}\label{fPhif}
 There exists an  equivalence
\begin{equation}\label{feqf}
%\begin{array}{ll} &
\fPhif:D^\LG(\N)\iso \Db(\fPf_I),
%;\\
%&\fPhif^{mix}:D^{\LG\times \Gm}(\N)\iso \Db(\fPfmI)
%\end{array}
\end{equation}
such that \begin{equation}\label{charactera}
%\begin{array}{ll}
pr_{ff}\circ \fPhi\cong  \fPhif\circ p_{Spr*},
 %pr_{ff}\circ \fPhi^{mix}\cong  \fPhif^{mix}\circ p_{Spr*}.
%\end{array}
\end{equation}
\end{Thm}

\begin{Rem}\label{derfu}
The functor $(\fPhif)^{-1}$ %, $(\fPhif^{mix})^{-1}$
is a derived functor of a left exact functor.

Namely, let $\unO_{\LG}$ be the ind-object of $Rep(\LG)$
corresponding to the module of regular functions on $\LG$, where
$\LG$ acts by left translations. %; here $H=\LG$ or $\LG\times \Gm$.
Notice that for  an object $\F$ of the derived category of $\LG$-equivariant
coherent sheaves one has $R^i\Gamma(\F)=Hom(\O, \F\otimes \unO[i])$.

For $X\in \fPfI$ %, $X'\in \fPfI^{mix}$ 
the space
$$Hom_{\fPfI}(\delta_e, X*\Z(\unO_\LG))=\oplusl V_\lambda^* \otimes
Hom(\delta_e ,X*Z_\lambda);$$
%$$Hom_{\fPfI^{mix}}(\delta_e, X*\Z(\unO_{\LG\times \Gm}))
%=\oplusl V_\lambda^* \otimes Hom(\delta_e,X*Z_\lambda(\frac{n}{2}))$$
can be given a structure of a $\LG$ %(respectively, $\LG\times \Gm$)
equivariant $\O(\N)$-module. Thus we get a left exact functor
$H:\fPfI\to Coh^\LG(\N)$; %, $\fPfI^{mix}\to Coh^{\LG\times \Gm}(\N)$
we claim that its derived functor $RH$ is isomorphic to
 $(\fPhif)^{-1}$. %, $(\fPhif^{mix})^{-1}$.

A sketch of the proof of this claim is as follows (the claim will not be used below, and details
of the proof are omitted). 
It follows from the Theorem and its proof that $\fPhif(V\otimes \O_\N)\cong \Z(V)\in \fPfI$, $V\in Rep(\LG)$,
and that $H\circ \fPhif |_{Coh_{fr}^\LG(\N)}\cong id_{Coh_{fr}^\LG(\N)}$ canonically, where
 $Coh_{fr}^{\LG}(\N)\subset Coh^\LG(\N)$ is the full subcategory consisting of objects
 of the form $V\otimes \O_\N$, $V\in Rep(\LG)$ (we will call such objects free sheaves). 
 
 The proof of  Theorem \ref{fPhif} below shows also that for  a finite complex $C^\bu$ of objects
of $Coh_{fr}^\LG(\N)$ the object $\fPhif(C^\bu)$ is represented by the complex
$(\fPhif(C^i))$.
 
Furthermore, it follows from Theorem 7 of \cite{1} that $Ext^{>0}_{\fPfI} (\delta_e, Z_\lambda)=0$, thus
for a complex $C^\bu$ of object of $\fPfI$ of the form $C^i=\Z(V_i)$, the object $RH(C^\bu)$
is represented by the complex $(H(C^i))$. Thus we get 
 a canonical isomorphism
\begin{equation}\label{RH}
RH\circ \fPhif |_{D_{fr}^\LG(\N) } \cong id | _{D_{fr}^\LG(\N)},
\end{equation} 
where $D_{fr}^\LG(\N) \subset D^\LG(\N)$  represented by a finite complex of free sheaves.

Finally, observe that the functor $RH\circ \fPhif$ sends $D^{<0}(Coh^\LG(\N))$ to itself, because:
\begin{multline*}\F\in D^{<0}(Coh^\LG(\N))\Rightarrow Ext^{\geq 0}(V\otimes \O, \F)=0 \ \forall V\in Rep(\LG)
\Rightarrow \\
 Ext^{>0} (\Z(V) ,\fPhif(\F))\ \forall V\in Rep(\LG) \Rightarrow RH(\fPhif(\F))\in D^{<0}
(Coh^\LG(\N)).\end{multline*}
Together with \eqref{RH} this shows that $RH\circ \fPhif \cong id$.

\end{Rem}

\begin{Rem}
Theorem \ref{fPhif} implies that the functor $Rp_*:D^\LG(\Nt)\to
D^\LG(\N)$ identifies $D^\LG(\N)$ with the quotient category
$D^\LG(\Nt)/Ker(Rp_*)$. 
%(and similarly for $Coh^\LG$ replaced by $Coh^{\LG\times \Gm}$). 

 Notice that the analogous statement with $D^b$ replaced by
$D^-$ is an immediate consequence of the isomorphism
$Rp_{Spr*}(\O_\Nt)\cong \O_\N$, and the fact that a triangulated
functor admitting a left adjoint which is also a right inverse is
 factorization by a thick subcategory, see e.g.  \cite{Verdier}, 
Proposition II.2.3.3.

It is  natural to ask whether the equivalence $D^\LG(\Nt)/Ker(Rp_*)
\cong D^\LG(\N)$ can be deduced directly from the isomorphism
 $Rp_*(\O_\Nt)\cong \O_\N$, and whether a similar equivalence holds
for an arbitrary proper morphism $p$ with $Rp_*(\O)\cong \O$. I do not know
the answer to these questions.
\end{Rem}

\subsection{Description of the $t$-structure on $D^\LG(\N)$}
 One can use the equivalences \eqref{feq}, \eqref{feqf} to
transport the tautological $t$-structure on the right-hand side to
a $t$-structure on the left-hand side. Let us call the resulting
$t$-structure on the derived category of equivariant coherent
sheaves the {\it exotic $t$-structure}. We provide
an explicit description of the exotic $t$-structure on
$D^\LG(\N)$. %, $D^{\LG\times \Gm}(\N)$.
\begin{Thm}\label{tstr}
The exotic $t$-structure on $D^\LG(\N)$ %, $D^{\LG\times \Gm}(\N)$
 coincides with the
perverse coherent $t$-structure corresponding to the
perversity given by
 \begin{equation}\label{pe}
p(O)=\frac{\codim(O)}{2};
\end{equation}
 see \cite{izvrat}, \cite{nilpokon}.
\end{Thm}

\subsubsection{}
Let $\orb$ denote the set of $\LG$-conjugacy classes of  pairs
$(N,\rho)$ where $N\in \N$, and $\rho$ is an irreducible
representation of the centralizer $Z_{\LG}(N)$.

 For a pair $(N,\rho)\in
\orb$ let $\L_\rho$ be the irreducible $\LG$-equivariant vector
bundle on the orbit $\LG(N)$, whose fiber at $N$ is isomorphic to
$\rho$. Let $j$ be the imbedding of $\LG(N)$ into $\N$. We have
the irreducible coherent perverse sheaf
$IC_{N,\rho}=j_{!*}(\L_\rho[\frac{-\codim \LG(N)}{2}])$, see
\cite{izvrat}.

\begin{Cor}\label{coodin}
a) We have $\fPhif(IC_{N,\rho})=L_w$ for some $w\in \fWf$.

b)  Identify $\Zet[W]$ with
the Grothendieck group $K(Coh^{\LG}(\St))$, where
$\St=\Nt\times_\N \Nt$ is the Steinberg variety, see \cite{CG},
\cite{Kthone}. Let $\Cbar_w$ be the Kazhdan-Lusztig basis of $\Zet[W]$
(specialization of the Kazhdan-Lusztig basis in the affine Hecke algebra at
$v=1$).
  Let $pr:\Zet[W]\to K(Coh^{\LG}(\N))$ be the
map induced by $Rp_*$, where $p:\St\to \N$ is the projection.

We have $pr(\Cbar_w)=0$ for $w\not \in \fWf$; and $pr(\Cbar_w)$ is the class
of an irreducible perverse coherent sheaf corresponding to the
perversity \eqref{pe}.

\end{Cor}

\begin{Rem} The Corollary implies validity of
 Conjectures 1 and 2, and first part of Conjecture 3 in \cite{Ostrik}.
\end{Rem}

\subsection{Dualities} We will denote the Verdier duality
functor on various categories  by $\bbV$. Thus $\bbV$ is an contravariant 
auto-equivalence of the abelian category $\PI$, which induces auto-equivalences
of the quotient categories $\fPI$, $\fPfI$ and their derived categories.

Define the anti-autoequivalence $\sigma$ %, $\sigma'$
 of $D^\LG(\N)$ %, $D^{\LG\times \Gm}(\N)$ respectively
by $\F\mapsto RHom(\F,\O)$.
 It is well-known that $\N$
is a Gorenstein scheme, and the dualizing sheaf for $\N$ is trivial.
 Thus $\sigma$ %, $\sigma'$ 
coincides with the Grothendieck-Serre
duality up to homological shift. 
%and (for $\sigma'$) a twist by a character of $\Gm$.
 Let $\kappa:\LG\to \LG$,
%$\kappa':\LG\times \Gm \to \LG\times \Gm$ 
be an automorphisms
which sends an element $g$ to an element conjugate to $g^{-1}$.
%(the Chevalley involution).
 We will also use the same letter to denote the induced push-forward functor
on the categories of representations and  equivariant coherent sheaves.

\begin{Thm}\label{dual}
We have \begin{equation}\label{duone}
 \fPhif\circ \kappa \circ
\sigma\cong \bbV \circ \fPhif;
\end{equation}
%\begin{equation}\label{dutwo}
 %\fPhif^{mix} \circ \kappa'
%\circ \sigma'\cong \bbV \circ \fPhif^{mix}.
%\end{equation}
\end{Thm}

%\subsubsection{} One can also describe the objects $(\fPhif^{mix})^{-1}(L_w)$.

\subsection{Cells and nilpotent orbits}\label{cel}
Recall the notion of a {\it two-sided cell} in $W$, and the
bijection between the set of two-sided cells in $W$
and nilpotent conjugacy classes in $\Lg$, see \cite{cells4};
for a two-sided cell $\c$ let $N_\c\in \Lg$ be a representative of the
corresponding conjugacy class.

For a two-sided cell $\c\subset W$ let
 $\P_I^{\leq \c}\subset \PI$ be the Serre subcategory generated by
irreducible objects  $L_w$, $w \in
\cupl_{\c'\leq \c} \c'$; and let $\fP_I^{\leq \c}\subset \fPI$,
 $\fPfI^{\leq \c}\subset \fPfI$
be the images of $\P_I^{ \leq \c}$. Let also $\Db_{\leq \c}(\fPI)
\subset \Db(\fPI)$,
$\Db_{\leq \c}(\fPfI)\subset \Db(\fPfI)$ be the full triangulated
 subcategories  generated by $\fP_I^{\leq \c}$,
$\fPfI^{\leq \c}$ respectively.
Replacing the nonstrict inequality by the strict one we get the definition
of categories %$\fP_I^{< \c}$,
 $\fPfI^{< \c}$, % $\Db_{< \c}(\fPfI)$,
$\Db_{< \c}(\fPfI)$ etc.

For a closed $\LG$-invariant subset $S\subset \N$ or $S\subset \Nt$
let $D_S^\LG(\N)\subset D^\LG(\N)$ (respectively, $D_S^\LG(\Nt)$)
be the full subcategory of complexes whose cohomology sheaves are
set-theoretically supported on $S$ (i.e. they are supported on some,
possibly nonreduced,
subscheme with topological space $S$). We abbreviate
 $D_{\leq N_\c}(\N)=D_{\overline{\LG(N_\c)}}^\LG(\N)$; 
 $D_{\leq N_\c}(\Nt)=D_{p_{Spr}^{-1}(\overline{\LG(N_\c)})}^\LG(\N)$.

% We will use the same notation for a functor on $\P$
%which preserves the above subcategories, and for the induced
%functor on $\P_\c$, $\A_\c$.

\begin{Thm}\label{kle} a) $\Db_{\leq \c}(\fP_I)
=\fPhi (D_{\leq N_\c}(\Nt))$;

$\Db_{\leq \c}(\fPfI ) =\fPhif(D_{\leq N_\c}(\N) )$.

b) We have $$c_1\leq c_2 \iff N_{\c_1}\in \overline{\LG(N_{\c_2})},$$
where the inequality in the left hand side refers to the standard
partial order on the set of 2-sided cells.
\end{Thm}

\begin{Rem} Part (b) of the Theorem was conjectured by Lusztig, see
 \cite{cells4}. %?
\end{Rem}

\subsection{Duflo involutions} Recall the notion of 
 a {\it Duflo} (or {\it distinguished})
{\it involution} in an affine Weyl group. We quote two of several available equivalent definitions.
On the one hand, an element $w\in W$ is a Duflo involution if and only if the corresponding element
in the asymptotic Hecke algebra (which is the Grothendieck ring of the truncated convolution
category, see the next subsection) is an idempotent. 
Moreover, the sum of all these idempotents over all Duflo involutions is the unit element in the
asymptotic Hecke algebra.

On the other hand, an element $w\in W$ is a Duflo
involution if and only if the degree of the Kazhdan-Lusztig polynomial $P_{e,w}$ is equal to $a(w)$,
where $a:W\to \Zet$ is Lusztig's $a$-function (recall that for any $w$ the degree
of $P_{e,w}$ is at most $a(w)$).
The latter characterization will be used in the proof
of Lemma \ref{Homnenol} below.

It is known that
for each two sided cell
$\c\subset W$ the set $\c\cap \fWf$ contains a unique Duflo involution, it will be denoted by
$d_\c$.

For a $\LG$-orbit $O\subset \N$ let $\Ohat_O$ denote
the sheaf $j_*(\O)$, where $j$ is the imbedding $O\imbed \N$,
and $j_*$ denotes the {\it non-derived} direct image.

\begin{Prop}\label{Duflo}
 We have $$\fPhif\left(\Ohat_{\LG(N_\c)}[-\frac{\codim \LG(N_\c)}{2}]
\right)\cong
L_{d_\c}.$$
\end{Prop}

\begin{Rem}
Proposition implies Conjecture 4
in \cite{Ostrik}.
\end{Rem}

\subsection{Truncated convolution categories}\label{trconv} 
In \cite{Lco} Lusztig defined for every two sided cell a monoidal category,
whose simple objects are $L_w$, $w\in \c$. He conjectured a relation
between this category and representations of the group $Z_\LG(N_\c)$;
these conjectures were partly proved in \cite{B}, \cite{BO}.
More precisely, one of the results of \cite{B} is as follows.
 Let $\PI^\c$ denote the Serre
quotient category $\PI^{\leq \c}/\PI^{< \c}$, and $\A_\c\subset
\PI^\c$ be the full subcategory consisting of subquotients of
objects of the form $\Z(V)*L_w \mod \PI^{<\c}$, $V\in Rep(\LG)$,
$w\in \c$. Let also $\A_\c^f\subset \A_\c$ be the subcategory 
 consisting of subquotients  $\Z(V)*L_w \mod \PI^{<\c}$, $V\in Rep(\LG)$,
$w\in \c\cap \fWf$.
Convolution with a central sheaf $\Z(V)$ induces
a functor on $\A_\c,$ $\A_\c^f$ which is also denoted by $X\mapsto \Z(V)*X$.

Truncated convolution provides $\A_\c$, $\A_\c^f$ with the structure
of a monoidal category. In \cite{B} we identified the monoidal category
$\A_\c^f$ with the
category of representations of a subgroup $H_\c\subset Z_\c$, where
 $Z_\c$
denotes the centralizer of $N_\c$ in $\LG$; in particular,
we have the restriction functor $r_\c^f:Rep(Z_\c)\to \A_\c^f$
(see Proposition \ref{ac} below for a more detailed statement).

We will compare $r_\c^f$ with a functor arising from $\fPhif$.
Set $\Db_\c(\fPfI):=\Db_{\leq \c}(\fPfI)/\Db_{<\c}(\fPfI)$; let
$Coh_{N_\c}^\LG(\N)$ be the category of equivariant coherent
sheaves on the formal neighborhood of the orbit $\LG(N_\c)$ in
$\N$, and $D_{N_\c}^\LG(\N)\cong D_{\leq
N_\c}^\LG(\N)/D_{<N_\c}(\N)$ be its bounded derived category.

 By Theorem \ref{kle}(a) the functor $\fPhif$ induces an
equivalence $D_{N_\c}^\LG(\N)\iso \Db_\c(\fPfI)$; we denote
this equivalence by $\Phi_\c$. 

\begin{Prop}\label{risPhi}
For $\rho\in Rep(Z_\c)$ we have a canonical isomorphism in $\fPfI^{\c}$
\begin{equation}\label{onorbit}
\Phi_\c(\L_\rho[-m])\cong r_\c^f(\rho),
\end{equation}
where $m=\frac{\codim \LG(N_\c)}{2}$.
\end{Prop}

\begin{Cor}\label{HZ}
We have $H_\c=Z_\c$.
\end{Cor}

\begin{Rem}
 A bijection between the set $\Lambda^+$ of dominant weights of $\LG$
 (which is the same as dominant coweights of $G$) and the set $\orb$
 was defined in \cite{nilpokon}; let $\iota_1$ denote this
bijection.
From the definition of $\iota_1$ in \cite{nilpokon}, it follows
that
%\begin{equation}\label{ione}
$$\fPhif(IC_{\iota_1(\lambda)})=L_{w_\lambda},
$$
% \end{equation}
where
$\{ w_\lambda\}=\fWf\cap W_f\cdot \lambda\cdot W_f.$

 Another bijection 
between the same sets (which we denote by $\iota_2$) was defined in
\cite{B}.
 $\iota_2$ is characterized as follows. 
If
$(N,\rho)=\iota_2(\lambda)$, and $N=N_\c$ for a two-sided cell
$\c$ then we have an isomorphism in $\PI^\c$
%\begin{equation}\label{itwo}
$$r_\c(\rho)\circ L_{d_\c}\cong L_{w_\lambda} .
$$
%\end{equation} 

Thus Proposition \ref{risPhi} implies that $\iota_1=\iota_2$.
\end{Rem}

\begin{Rem}
The equality $\iota_1=\iota_2$ implies  Conjecture 3 in \cite{Ostrik}.
\end{Rem}

\end{section}

\section{Proofs}
\subsection{Proof of Theorem \ref{fPhif}}
\begin{Lem}\label{fPf}
For $w\not \in \fWf$
 we have $p_{Spr*}(\fPhi^{-1}(L_w))=0$.
\end{Lem}

\proof For $\F\in D^\LG(\Nt)$ we have $p_{Spr*}(\F)=0$ iff
$Ext^\bu (V_\lambda\otimes \O, \F)=0$ for all $\lambda\in
\Lambda^+$. Thus we need to check that for $X\in \Db(\fPI)$ we
have
\begin{equation}\label{ZLw}
Ext^\bu_{\fP}(Z_\lambda, L_w)=0
\end{equation}
for $w\not \in \fWf$.

We will check the equivalent statement
%\begin{equation}\label{iff}
$$
Ext^\bu_{\DIW}(\Delta_e*Z_\lambda,\Delta_e*L_w)=0
$$
%\end{equation}
for $w\not \in \fWf$;  cf \cite{1}, Theorem 2.

If $w \in \fW$ but $w\not \in \fWf$ then for some simple root
$\alpha$,  $\alpha\ne
\alpha_0$ we have $L_w=\pi_\alpha^!(L'_w)$;
here $\alpha_0$ is the affine simple root, $\pi_\alpha:\Fl\to \Fl(\alpha)$
is the projection ($\Pone$ fibration) to the corresponding
partial affine flag variety, and $L_w'$
is an $\bI$-equivariant constructible complex on $\Fl(\alpha)$
(actually, $L'_w[1]$ is a perverse sheaf).
Then $\Delta_e*L_w\cong \pi_\alpha^!(\Delta_e*L'_w)$, and
we have $$Ext^\bu(\Delta_e*Z_\lambda,\Delta_e* L_w)=
Ext^\bu({\pi_{\alpha!}}(\Delta_e*Z_\lambda), \Delta_e*L_w')=0,$$
because ${\pi_{\alpha!}}(\Delta_e*Z_\lambda)=0$, cf. \cite{1}, proof of Lemma
28. \epf

\subsubsection{}
The Lemma shows that the functor $p_{Spr*} \circ
\fPhi^{-1}:\Db(\fPI) \to D^\LG(\N)$
%;  $p_{Spr*} \circ (\fPhi^{mix})^{-1}:\Db(\fPI) \to D^\LG(\N)$ 
factors through
$\Db(\fPfI)$ %, $\Db(\fPfI^{mix})$ respectively 
(here we use that
$\Db(\A/\B)\cong \Db(\A)/\Db_\B(\A)$ for an abelian category $\A$,
and a Serre subcategory $\B$, where $\Db_\B(\A)\subset \Db(\A)$ is
the full subcategory of objects with cohomology in $\B$).
%, the quotient
%in the left hand side is the Serre quotient, and the quotient in the
%right-hand side is the quotient of a triangulated subcategory by a thick
%subcategory). see ??.

It remains to check that the resulting functor
$$\Upsilon: \Db(\fPfI)\to D^\LG(\N),$$
%$$\Upsilon_{mix}: \Db(\fPI^{mix})\to D^{\LG\times \Gm}(\N)$$
is an equivalence; then $\fPhif:=\Upsilon^{-1}$
%,$\fPhif^{mix}:=\Upsilon^{-1}_{mix}$
 clearly satisfies the conditions
of the Theorem. %We will argue that $\Upsilon$ is an equivalence;
%the proof for $\Upsilon_{mix}$ is parallel.

\subsubsection{}
Let us  check that $\Upsilon$ is a full imbedding.

First we claim that
\begin{equation}\label{ism}
Hom_{\Db(\fPfI)}(X,Y)\overset{\Upsilon}{\To} Hom_{D^\LG(\N)}
(\Upsilon(X),\Upsilon(Y))
\end{equation}
is an isomorphism for $X=Z_\lambda$. Indeed, \eqref{ZLw} implies that
$Hom_{\Db(\fPI)}(Z_\lambda,Y)\iso Hom_{\Db(\fPfI)}(Z_\lambda,Y)$.
Also, the equality $Rp_{Spr*}(\O_{\Nt})=\O_\N$ implies that
$$\begin{array}{ll}
Hom_{D^\LG(\Nt)}(V_\lambda\otimes \O_\Nt,\F)\iso
 Hom_{D^\LG(\Nt)
}(Rp_{Spr*}(V_\lambda\otimes \O_\Nt),Rp_{Spr*}(\F))=
\\
 Hom_{D^\LG(\N)}(V_\lambda\otimes \O_\N,Rp_{Spr*}(\F))
\end{array}
$$
for $\F\in D^\LG(\Nt)$.
Thus validity of \eqref{ism} for $X=Z_\lambda$ follows from
$\fPhi^{-1}$ being an equivalence.

We now want to deduce that \eqref{ism} is an isomorphism for all $X$.
%$X=j_{w*}$.
 The argument is a version of the  proof of the fact
that an effaceable $\delta$-functor is universal.

\begin{Lem}\label{5}
Let $D$ be a triangulated category, $F=(F^i)$, $F'=({F'}^i)$ be cohomological
functors from $D$ to an abelian category, and $\phi:F\to {F'}$ be
a natural transformation. Let $S\subset D$ be a set
of objects. Assume that

i) There exists $d\in \Zet$ such that $F^i(X)=0={F'}^i(X)$ for $i<d$,
$X\in S$.

ii) For any $X\in S$ there exists an exact triangle
$X\to \Xtil \to Y$
where $Y\in S$, and $\phi: F^i(\Xtil)\to {F'}^i(\Xtil)$
is an isomorphism for all $i$.

Then $\phi:F^i(X)\to {F'}^i(X)$ is an isomorphism for all $X\in S$.
\end{Lem}

\proof We go by induction in $i$. Condition
 (i) provides the base of induction. Applying the 5-lemma to
$$\CD
F^{i-1}(\Xtil)@>>> F^{i-1}(Y)@>>> F^i(X) @>>> F^i(\Xtil) @>>> F^i(Y) \\
@V{\|}VV  @V{\|}VV @VVV @V{\|}VV @VVV \\
{F'}^{i-1}(\Xtil)@>>> {F'}^{i-1}(Y)@>>> {F'}^i(X) @>>>
 {F'}^i(\Xtil) @>>> {F'}^i(Y)
\endCD
$$
we see that $\phi:F^i(X)\imbed {F'}^i(X)$. Since $Y\in S$ we have
$\phi:F^i(Y)\imbed {F'}^i(Y)$ which implies $\phi:F^i(X)\iso {F'}^i(X)$. \epf

To exhibit a generating  set for $\Db(\fPfI)$
satisfying the conditions of Lemma \ref{5}
we need another Lemma.

A filtration on the object of $\fPI$ (respectively, $\fPfI$) will be
called {\it costandard} if its associated graded is a sum of objects
$j_{w*}$, $w\in \fW$ (respectively, $w\in \fWf$).
Such a filtration will be called {\it standard}
 if its associated graded is a sum of objects
$j_{w!}$, $w\in \fW$ (respectively, $w\in \fWf$).

\begin{Lem}\label{YZX}
a) If $w_1, w_2\in W$ and $w_2\in W_f\cdot w_1\cdot W_f$, then
$j_{w_1*}$ and $j_{w_2*}$ are isomorphic in $\fPfI$.
% ; also $j_{w_1!}$ and $j_{w_2!}$ are isomorphic in $\fPfI$.

b) Let $X\in \fPfI$ be an object with a costandard filtration.
Then there exists a short exact sequence $0\to Y\to Z \to X \to 0$
in $\fPfI$ where $Z$ is a (finite) sum of objects $Z_\lambda$,
 and $Y$ has a costandard filtration.
\end{Lem}

\proof (a) %We check the first statement; the second one then follows
%by duality.
%
We can assume that $w_2=s w_1$ or $w_2=w_1s$
for a simple reflection $s=s_\alpha\in W_f$, and that $\ell(w_2)>\ell(w_1)$.
Assume first that $w_2=sw_1$.
We have $j_{w_2*}=j_{s*}*j_{w_1*}$. The short exact sequence
$$0\to \delta_e \to j_{s*}\to L_s\to 0$$
(where $e\in W$ is the identity, and $\delta_e=j_{e*}=j_{e!}=L_e$
is the unit object of the monoidal category $D_I(\Fl)$)
yields an exact triangle
$$j_{w_1*}\to j_{w_2*}\to L_s *j_{w_1*}.$$
 It is easy to see
that   $L_s*j_{w_1*}$ is a perverse sheaf;
this object is equivariant with respect to
the parahoric group scheme $\bI_\alpha$. It follows that any its irreducible
subquotient is also  equivariant under
this group; hence such a subquotient
is isomorphic to $L_w$ for some $w$ satisfying
 $\ell(sw)<\ell(w)$. This shows that
$L_s*j_{w_1*}$ is zero in $\fPfI$, hence
$j_{w_2*}$ and $j_{w_1*}$ are isomorphic in $\fPfI$.

In the case $w_2=w_1s$ the proof is parallel, with the words
``is equivariant under $\bI_\alpha$'' replaced by ``lies in the
image of the functor  $\pi_\alpha^*$''. Thus (a) is proved.

(b) Theorem 7 of \cite{1} implies that $Z_\lambda$ (considered
as an object of $\fPI$) has both a standard and a costandard filtration.
The top of the latter is a surjection $f_\lambda:
Z_\lambda\to j_{\lambda*}$
whose kernel, subsequently, has a costandard filtration.
Taking the image of $f_\lambda$ in $\fPfI$ we get an arrow
$\fbar_\lambda:Z_\lambda\to j_{w_\lambda*} $ in $\fPfI$
whose kernel has a costandard filtration by (a); recall that $\{w_\lambda\}
=\fWf\cap W_f\cdot \lambda \cdot W_f$.

Let now $X\in \fPfI$
be an object with a costandard filtration; let $0\to X' \to X\to j_{w*}\to 0$
be the top of the filtration. By induction in the length of the filtration
we can assume the existence of an exact sequence
$$0\to Y' \to Z' \overset{f'}{\To} X'\to 0$$ of the required form.
We have $w=w_\lambda$ for some $\lambda\in \Lambda^+$. We claim that the
surjection $\fbar_\lambda:Z_\lambda\to j_{w*}$ factors through a map
$Z_\lambda\to X$. Indeed, the obstruction lies in $Ext^1_{\fPfI}(Z_\lambda,
X').$ We claim that
$$Ext^1_{\fPfI}(Z_\lambda, j_{w*})= Ext^1_{\fPI}(Z_\lambda, j_{w*})=0.$$
Here the first equality follows from  \eqref{ZLw}. The second one
is a consequence of the existence of a standard filtration
on $Z_\lambda$ (considered as an object of $\fPI$), and the equality
$$Ext_{\fPI}^\bu(j_{w_1!}, j_{w_2*})=\Ql^{\delta_{w_1,w_2}}.$$
The latter is a consequence of \cite{1}, Theorem 2 and Lemma 1,
which identify the left-hand side with an Ext space in the derived
category of $l$-adic complexes on $\Fl$ (more precisely, with
$Ext^\bu(\Delta_{w_1}, \nabla_{w_2})$ in the notations of \cite{1}).

Now let $\tilde f_\lambda:
 Z_\lambda\to X$ be some map, such that the composition
$Z_\lambda\to X\to j_{w*}$ equals $\fbar_\lambda$.
Then we set $Z=Z'\oplus Z_\lambda$, and the map $f:Z\to X$ is set
to be $f:=f'\oplus \tilde f_\lambda$. The exact sequence
$$0\to Ker(f')\to Ker(f)\to Ker(\fbar_\lambda)\to 0$$
shows that $f$ satisfies the requirements of (b).
\epf

\subsubsection{} We can now finish the proof of $\Upsilon$
being an equivalence.
 We apply Lemma \ref{5} to the following data:

\noindent
 $D=\Db(\fPfI)^{op}$,

\noindent
$F:X\mapsto Hom^\bu(X,X_0)$ for some fixed $X_0\in \Db(\fPfI)$;

\noindent
$F':X\mapsto Hom^\bu(\Upsilon(X), \Upsilon(X_0))$;

\noindent
 the transformation $\phi$
comes from functoriality of $\Upsilon$;

\noindent
the set $S$ consists of all objects
of $\fPfI$ which have a costandard filtration.

We claim that conditions of Lemma \ref{5} are satisfied. 

In fact, condition (i) is satisfied for any $d$ such that $X_0\in D^{\geq -d}(\fPfI)$,
$\Upsilon(X_0)\in D^{\geq -d}(Coh^\LG(\N))$. Vanishing of $F^i(X)$, $X\in S$, $i<d$ follows then
from vanishing of negative Ext's between objects of $\fPfI$, while vanishing for $(F')^i(X)$,
   $X\in S$, $i<d$ follows from vanishing of negative Ext's in $Coh^\LG(\N)$, in view of the fact that
 $\Upsilon (S) \subset Coh^\LG(\N)$. The latter inclusion amounts to the fact that
 $Ext^{>0}_{\fPI}(\Z(V), X)=0$, $X\in S$, $V\in Rep(\LG)$, 
 which follows from the tilting property of central sheaves, Theorem 7 of \cite{1}.

Since \eqref{ism} has been proven to be an isomorphism for $X=Z_\lambda$, 
Lemma
\ref{YZX}(b) shows that condition (ii) of Lemma \ref{5} is satisfied.
Hence \eqref{ism} is an isomorphism whenever $X$ has a costandard
filtration; in particular, for $X=j_{w*}$. But $\Db(\fPfI)$ is
generated as a triangulated category by $j_{w*}$, $w\in \fWf$;
hence \eqref{ism} holds for all $X$, i.e. $\Upsilon$ is a full
imbedding.

It remains to show that $\Upsilon$ is essentially surjective;
since it is a full imbedding it suffices to see that the image of
$\Upsilon$ contains a set of objects generating $D^\LG(\N)$ as a
triangulated category. This is done in Lemma 7 of \cite{nilpokon}.
\epf

\subsection{Proof of Theorem \ref{tstr}}
It follows from the results of  \cite{nilpokon} that $D^\LG(\N)$ carries a unique
$t$-structure such that the objects $A_\lambda$ lie in
the heart of this $t$-structure for all $\lambda$, where
$A_\lambda=p_{Spr*}(\O(\lambda))$; and this $t$-structure
coincides
 with the perverse
coherent $t$-structure corresponding to the perversity function $p(O)=\frac{\codim O}{2}$
(which coincides with the middle perversity up to a total shift by $\frac{\dim \N}{2}$).
Recall that the objects  $J_\lambda\in \PI$ (the Wakimoto sheaves),
satisfy $\fPhi(\O(\lambda))\cong J_\lambda$; hence
\begin{equation}\label{kudyA}
\fPhif(A_\lambda)\cong pr_{ff}(J_\lambda).
\end{equation}
Thus the heart of the $t$-structure obtained by transport of the tautological $t$-structure on $D^b(\fPfI)$ under the equivalence $(\fPhif)^{-1}$ contains the objects $A_\lambda$ in its heart, so it coincides
with the perverse coherent $t$-structure. \epf

%This proves the statement concerning $D^\LG(\N)$. The one
%concerning $D^{ \LG\times \Gm}(\N)$ follows, because
%$\fPhif^{mix}$ intertwines the natural functor $\fPfI^{mix}\to
%\fPfI$ with the forgetful functor $D^{\LG\times \Gm}(\N)\to
%D\LG(\N)$; and the perverse $t$-structure on $D^{\LG\times
%\Gm}(\N)$ is uniquely characterized by the fact that the forgetful
%functor sends it to the perverse coherent $t$-structure on
%$D^\LG(\N)$.

\subsubsection{Proof of Corollary \ref{coodin}}
a) is immediate from Theorem \ref{tstr}, because $IC_{N,\rho}$ is
an irreducible object in the heart of the perverse coherent
$t$-structure, so $\fPhif(IC_{N,\rho})$ is an irreducible object
of $\fPfI$. Let us prove (b). Let $\gp$ denote the map
on Grothendieck groups induced by the composition
$$\Db(\PI)\to \Db(\fPfI)\overset{\fPhif^{-1}}{\To} D^\LG(\N).$$
We can identify $\Zet[W]$ with $K(\Db(\PI))$ by means of the isomorphism
 sending $(-1)^{\ell(w)}
\cdot w$ to the class  $[j_{w!}]=[j_{w*}]$;
it maps $\Cbar_w$ to the class of $L_w$.
Thus it is clear that $\gp(C_w)=0$ for $w\not \in \fWf$; and (a) shows that
$\gp(C_w)$ is the class of an irreducible perverse coherent sheaf.
It remains to check that $\gp=pr$.
This follows from: $\gp(w)=(-1)^{\ell(w)}[A_\lambda]=pr(w)$ for $\lambda\in \Lambda^+$,
$w\in W_f\cdot \lambda \cdot W_f$.  Here the first equality follows from 
 \eqref{kudyA} for $w=\lambda$ and from Lemma \ref{YZX}(a) in general. The second equality
 holds by
  Lemma 2.4 in \cite{Ostrik}.
  \epf

%We have
%$$pr_{ff}(j_{w*}) %\cong pr_{ff}(j_{\lambda_w*})\cong 
%A_\lambda,$$
%where the first isomorphism follows from  Lemma \ref{YZX}(a),
%and the second one is \eqref{kudyA} above.
%Thus $\gp(w)=(-1)^{-\ell(w)+\ell(\lambda_w)}[A_\lambda]=pr(w)$.

\subsection{Proof of Theorem \ref{dual}} 
%We write down
%the argument for the ``unmixed'' case; the ``mixed'' case is
%similar.
%
%
Recall that $\Sa$ denotes the equivalence $Rep(\LG)\to
\P_\GO(\Gr)$. Let $\upsilon:Rep(\LG)\to Rep(\LG)^{op}$ denote the
functor $V\mapsto V^*$.

 Recall also that $Coh_{fr}^\LG(\N)$
denotes the category of $\LG$-equivariant vector bundles on $\N$
which have the form $V\otimes \O$, $V\in Rep(\LG)$. Thus
$Coh_{fr}^\LG(\N)$ is a tensor category under the tensor product
of vector bundles. It was shown in \cite{1} that the map $V\otimes
\O\mapsto \Z(V)$ extends naturally to a monoidal functor
$Coh_{fr}^\LG(\N)\to D_I(\Fl)$; we denote the resulting monoidal
functor by  $\tilde \Z$
(thus $\tilde \Z=F\circ p_{Spr}^*$ in notations of \cite{1}).

\begin{Lem}\label{VeZ} 
%a) There exists a tensor isomorphism of functors
%$Rep(\LG)^{op}\to \P_\GO(\Gr)$
%$$\Sa \circ ( \upsilon \circ \kappa) \cong \bbV \circ \Sa.$$
%
%b)
 There exists a tensor isomorphism of functors $Coh_{fr}^\LG(\N)^{op}\to
\P_\GO(\Gr)$ $$ \tilde\Z \circ (\sigma \circ \kappa ) \cong \bbV
\circ \tilde \Z.$$

\end{Lem}

\proof 
%a) Any tensor functor $F:Rep(\LG)\to Rep(\LG)$ such that
%$F(V)\cong V$ for any $V\in Rep(\LG)$ is isomorphic to identity
%(indeed, by Tannakian formalism it is induced by some automorphism
%of the group $\LG$; but it is well-known that an automorphism of a
%connected reductive group over an algebraically closed field,
%inducing the identity map on the set of isomorphism classes of
%irreducible representations is inner). Applying this observation
%to the functor $\Sa^{-1}\circ \bbV \circ \Sa \circ (
%\upsilon\circ \kappa)$ we get the statement.
%
%b) 
The functor $\tilde \Z$ is characterized  by the following two
conditions (cf. \cite{1}, Proposition 4(a)):
\begin{equation}\label{Zone}
\tilde\Z|_{Rep(\LG)}\cong \Z
\end{equation}
\begin{equation}\label{Ztwo}
\tilde\Z(N^{taut}_{V\otimes \O})=\M_{\Z(V)}.
\end{equation}
More precisely, given a functor $\tilde\Z '$ with a functorial tensor
isomorphism $\tilde \Z'(V\otimes \O)\cong \Z(V)$, $V\in Rep(\LG)$,
 which intertwines
$\tilde\Z'(N^{taut}_{V\otimes \O})$ and $\M_{\Z(V)}$, one can construct
a canonical isomorphism $\tilde \Z'\cong \tilde \Z$.
Here $N^{taut}_\F$ is the "tautological" endomorphism of an
equivariant sheaf $\F\in Coh^\LG(\N)$, whose action on the fiber
at a point $x\in \N$ coincides with the action of $x\in
Stab_{\Lg}(x)$ coming from the equivariant structure; and
$\M_{\Z(V)}$ is the logarithm of monodromy endomorphism of $\Z(V)$
(arising from the construction of $\Z(V)$ via the nearby cycles
functor).
Thus we will be done if we show that \eqref{Zone} can be constructed, so that
 \eqref{Ztwo} holds,
 for $\tilde \Z$ replaced by $\bbV\circ \tilde \Z \circ
(\sigma \circ \kappa)$.

This follows from  
%(a) and the fact that nearby cycles commute with
%Verdier duality, by inspection of the definition of $\Z$ in
%\cite{KGB}. Equality \eqref{Ztwo} follows from 
existence
of natural isomorphisms satisfying the corresponding equalities:
\begin{equation}\label{ka}
\kappa(V\otimes \O)\cong \kappa(V)\otimes \O,\ \ \ \ \ \
\kappa(N^{taut}_{V\otimes \O})=N^{taut}_{\kappa(V)\otimes \O};
\end{equation}
\begin{equation}\label{si}
\sigma(V\otimes \O)\cong V^*\otimes \O,\ \ \ \ \ \
\sigma(N^{taut}_{V\otimes \O})=N^{taut}_{V^*\otimes \O};
\end{equation}
\begin{equation}\label{bbv}
\bbV(\Z(V))\cong \Z(V), \ \ \ \ \ \bbV(\M_{\Z(V)})=\M_{\Z(V)}.
\end{equation}
Here \eqref{ka} and \eqref{si} is an easy exercise; and
\eqref{bbv} follows from the fact nearby cycles commute with
Verdier duality, and the isomorphism $\bbV \circ \Psi\cong
\Psi\circ \bbV$ (where $\Psi$ is the nearby cycles functor)
respects the monodromy action, by inspection of the definition of $\Z$ in
\cite{KGB}. \epf

\subsubsection{} We are now ready for the proof of the Theorem.
Using the monoidal functor $\tilde \Z$ one can define the action
of the monoidal category $Coh_{fr}^\LG(\N)$ on $\fPfI$ and on its
derived category. It follows from the isomorphisms \eqref{charactera},
and the definition of $\fPhif$ in \cite{1} (cf. \cite{1}, Theorem 1)
 that $\fPhif$ intertwines this action with the
action of $Coh_{fr}^\LG(\N)$ on $D^\LG(\N)$ by tensor products;
i.e. we have a natural isomorphism
\begin{equation}\label{PhiZ}\fPhif(V\otimes \F)\cong \Z(V) *\fPhif(\F).
\end{equation}

 Set $\phi= 
% \kappa\circ \sigma \circ \fPhif^{-1} \circ bbV$.
\fPhif^{-1}\circ \bbV \circ \fPhif \circ \sigma\circ \kappa$. 
We want to construct an isomorphism $\phi\cong id$.

  Lemma \ref{VeZ} shows that $\phi$ commutes with the action of
$Coh_{fr}^\LG(\N)$ by tensor products. Furthermore, it is easy to
see that $\phi(\O)\cong  \O$. Thus we get an isomorphism
\begin{equation}\label{onfr}
\phi|_{Coh_{fr}^\LG(\N)}\cong id .\end{equation} 

Let now $C^\bu$ be a bounded complex where $C^i\in Coh_{fr}^\LG(\N)$,
and $C$ be the corresponding object of $D^\LG(\N)$.
By inspection of the definition of 
$\fPhif$ one checks that $\phi(C)$ is represented by the complex $(\phi(C^i))$.
%This follows by
%inspection of the definition of $\fPhif$, together with
%the fact  that
%Verdier duality applied to a constructible complex represented
%as a finite complex of perverse sheaves can be computed by applying
%the duality to the complex termwise \cite{Beil}.
This yields an isomorphism
\begin{equation}\label{onperf}
\phi|_{D^\LG_{fr}(\N)}\cong id
\end{equation} (see Remark \ref{derfu} for notation).
As in Remark \ref{derfu} we see that $\phi$ preserves
$D^{<0}(Coh^\LG(\N))$, which, together with \eqref{onperf}, yields
an isomorphism $\phi\cong id$.
 \epf

\subsection{Proof of Theorem \ref{kle}}\label{24}
We first recall some results of \cite{BO}, \cite{B}
(see section \ref{trconv} above for notations).

\begin{Prop}\label{ac} $\A_\c$ carries a natural structure of a
rigid monoidal category (given by the truncated convolution
$\circ$); $\A_\c^f\subset
\A_\c$ is a monoidal subcategory.
 Let $\bbone$  be the unit object\footnote{It follows from the results of Lusztig (cf. \cite{B})
 that $\bbone \cong L_{d_\c}$; Proposition \ref{Duflo} provides a description
 of the corresponding object in the derived category of coherent sheaves.} of $\A_\c$.
 We have a monoidal
central\footnote{See e.g. \cite{1}, \S 3.2, or (with more details) 
\cite{B}, \S 2.1 for a definition.}
 functor $r_\c:Rep(Z_\c)\to \A_\c$ such that

i) The composition of the restriction functor
$res^\LG_{Z_\c}:Rep(\LG) \to Rep(Z_\c)$ with $r_\c$ is isomorphic
to the functor $V \mapsto \Z(V)* \bbone$.

ii) The element $N_\c$ yields a tensor endomorphism of the functor
$Res^\LG_{Z_\c}$. The isomorphism of (i) carries this endomorphism
into the endomorphism induced by the logarithm of monodoromy, see
\cite{KGB}, Theorem 2.

iii) For $X\in \A_\c$ we have a canonical isomorphism
\begin{equation}\label{rc}
Z_\lambda*X\cong r_\c( V_\lambda|_{Z_\c}) \circ X.
\end{equation}

iv) The functor  $V\mapsto r_\c(V)\circ X$ from $Rep(Z_\c)$ to
$\A_\c$ is exact and faithful for all $X\in \A_\c$, $X\ne 0$.

v) The functor $r_\c^f$ defined by $r_\c^f(X)=r_\c(X)\circ L_{d_\c}$
is a monoidal functor $Rep(Z_\c)\to \A_\c^f$.

There exists an algebraic subgroup $H_\c\subset Z_\c$, and an equivalence
$\A_\c^f\cong Rep(H_\c)$, which intertwines $r_\c^f$ with the
restriction functor $Rep(Z_\c)\to Rep(H_\c)$.
\end{Prop}
\proof See \cite{BO}. \epf
%%%%

We need to spell out compatibility between \eqref{rc} and 
equivalence $\fPhi$.

The functor $F$ induces a map $$Hom_{Coh^\LG(\N)}(V_1\otimes \O,
V_2 \otimes \O) \to Hom (\Z(V_1), \Z(V_2));$$ where $V_1,V_2\in Rep(\LG)$.
For $h\in
Hom_{Coh^\LG(\N)}(V_1\otimes \O, V_2 \otimes \O)$ define
$h_X:\Z(V_1)*X \to \Z(V_2)*X$ by $h_X=F(h)*id_X$.

On the other hand, given $h\in Hom_{Coh^\LG(\N)}(V_1\otimes \O,
 V_2 \otimes \O)$ we can consider
the induced map of fibers at $N_\c$; we denote this map by
$h_{N_\c}\in   Hom_{Z_\c}(V_1, V_2)$.

\begin{Lem}\label{compaN}
Let $X\in \PI^{\leq \c}$, $X\mod \PI^{<\c}\in \A_\c$. Then for
 $h\in  Hom_{Coh^\LG(\N)}(V_1\otimes \O, V_2 \otimes
\O)$  isomorphism \eqref{rc} carries $h_X$ into
$r_\c(h_{N_\c})\circ id_X$.
\end{Lem}

\proof 
We need to enhance \eqref{rc}
to an isomorphism between the two actions of the tensor
category $Coh_{fr}^\LG(\N)$ on $\A_\c$, where the first one is given by
$\F:X\mapsto \tilde \Z(\F)*X$, while the second one is given by
$\F: X\mapsto r_\c(\F_{N_\c})\circ X$, where $\F_{N_c}$ denotes
the fiber of $\F$ at $N_\c$.
We apply   the (easy) uniqueness part of the
Proposition 4(a) in \cite{1} to the situation where the target category
$\C$ is the category of endo-functors of $\A_\c$. According to that 
Proposition, it suffices to check that \eqref{rc} is compatible
with the image of the tautological
endomorphism  $N^{taut}$ of $id_{Coh_{fr}^\LG(\N)}$. 
In view of Proposition \ref{ac}(iii), this compatibility
follows by comparing
Proposition \ref{ac}(ii) with  compatibility \eqref{Ztwo}
between $N^{taut}$ and monodromy via $\tilde \Z$. \epf

Theorem \ref{kle} will be deduced from the next
\begin{Lem} a)
  For $w\in  \fW$ we have
\begin{equation}\label{suppLw}
 w\in \c \ \ \Rightarrow \ \  p_{Spr} (supp(\fPhi^{-1}(L_w)))= \overline{\LG
 (N_\c)}.
\end{equation}

 b) For any
$X\in \Db(\fPI)$ we have
\begin{equation}\label{suppX}
  p_{Spr} (supp(\fPhi^{-1}(X)))= \cupl_{\c}\overline{\LG
 (N_\c)}
\end{equation}
where $\c$ runs over the set of such 2-sided cells that the
multiplicity of $L_w$ in the Jordan-Hoelder series of $H^i(X)$ is
non-zero for some $w\in \c\cap \fW$.
\end{Lem}
\proof Let $\cJ \subset \O_\N$ be the ideal sheaf of the closure
of a $\LG$-orbit $O$ on $\N$. Fix $n>0$.
 There exists a surjection of equivariant  sheaves
$V\otimes \O \sur \cJ^n $ for some $V\in Rep(\LG)$. Let $\phi:
V\otimes \O \to \O$ be the composition $V\otimes \O \sur \cJ^n
\imbed \O$; we use the same symbol to denote the pull-back of $\phi$
under $p_{Spr}$.
Then an object $\F\in D^\LG(\Nt)$ lies in $D_{p_{Spr}^{-1}( \Obar)}^\LG(\Nt)$ if
and only if the arrow $\phi\otimes id_\F: V\otimes \F \to \F$
equals zero for some (equivalently, for all large) $n$. Thus to
check \eqref{suppLw} it is enough to show that for $w\in \c$ we
have
\begin{equation}\label{va}
\overline{O}\owns N_\c \iff 0=\phi_{L_w}\in Hom(
\Z(V)*L_w , L_w). \end{equation}

If $w\in \c\cap \fW$ then a morphism $\Z(V)*L_w \to L_w$ is zero iff the
induced arrow in $\A_\c$ is zero; this is also equivalent
to the induced arrow in $\fPI$ being zero. In view of Lemma \ref{compaN}
the induced map $(\phi)_{L_w}\mod \PI^{<\c}\in Hom (\Z(V)*L_w ,
L_w) $ equals $r_\c(\phi_{N_\c})\circ id_{L_w}$. But
$\phi_{N_\c}=0$ if $N_\c\in \overline{O}$, so \eqref{va} holds in
this case. Conversely, if $\overline{O}\not \owns N_\c$ then $
\phi_{N_\c}\ne 0$ for all $n$. Since the functor $V\mapsto
r_\c(V)\circ X$ from $Rep(Z_\c)$ to $\A_\c$ is exact and faithful
for all $X\in \A_\c$, $X\ne 0$ we see that $\phi_{L_w}$
 is nonzero in this case. This shows \eqref{va}, and
hence \eqref{suppLw}.

\eqref{suppLw} implies that the left hand side of \eqref{suppX} is
contained in the right-hand side. Let us check the other
inclusion. Let $\J$ be the ideal sheaf of a proper $\LG$-invariant
subvariety $S$ in the right-hand side of \eqref{suppX}, and $\phi
:V\otimes \O\to \O$ satisfy $im(\phi)=\J^n$ as before.
 We need to verify that $supp(X)\not \subset S$, which is equivalent
to saying that
 the induced morphism $\Z(V)(X) \to X$ is
 nonzero. There exists $w\in \c\subset W$ such that  the
multiplicity of $L_w$ in the Jordan-Hoelder series of $H^i(X)$ is
non-zero for some $i$ but $N_\c\not \in S$. We saw in the previous
paragraph that the  morphism $(\phi)_{L_w}:\Z(V)*L_w\to L_w$ is
non-zero. But the latter is a subquotient of $H^i((\phi)_X)$; so
$(\phi)_X\ne 0$ as well. \epf

\subsubsection{Proof of Theorem \ref{kle} (conclusion)}
(a) follows from (b) and
\eqref{suppX}; so let us prove (b). Let $\c_1, \c_2\subset W$ be
two sided cells. Let $\J_i\subset \O_\N$ be the ideal sheaf of
$\overline{\LG(N_{\c_i})}$, and
 $\phi_i:V_i\otimes \O_\N\to \O_\N$ have $\J_i$ as its image ($i=1,2$).

Assume that $c_1\leq \c_2$; pick $w_1\in \c_1\cap \fW$, $w_2\in \c_2
\cap \fW$.
Then $L_{w_1}$ is a direct summand in the convolution
$X_1*L_{w_2}*X_2$ for some semisimple complexes $X_1, X_2 \in
\Db_I(\Fl)$. Hence the arrow $(\phi_2)_{L_{w_1}}$ is a direct
summand in $$X_1*((\phi_2)_{L_{w_2}}) *X_2= (\phi_2)_{X_1*L_{w_2}
*X_2}.$$  % where the equality follows from \eqref{compasve}.

%It is easy to check that for $X\in \Db(\P)$, $Y\in \Db_I(\P)$ we
%have \begin{equation}\label{compasve}
%h_{X*Y}=h_X*id_Y=id_X*h_Y.
%\end{equation}

But
$$(\phi_2)_{L_{w_2}}=0;$$
hence
$$(\phi_2)_{L_{w_1}}=0,$$
which implies $$N_{\c_1}\in p_{Spr}(supp(\fPhi ^{-1}(L_{w_1})) )\subset
\overline{\LG(N_{\c_2})}.$$
 Conversely, suppose that $N_{\c_1}\in
 \overline{\LG(N_{c_2})}.$ Let $$\bbK=(0\to
\Lambda^d(V)\otimes \O_\N \to \cdots \to V\otimes \O_\N\to \O_\N
\to 0)$$ be the Koszul complex of $\phi_1$. Pick $w\in \c_2\cap^fW$. Then
we have $$\overline{\LG(N_{\c_1})} = p_{Spr}(supp (\bbK \otimes _{\O_\N}
\fPhi^{-1}(L_w))).$$ Hence, according to \eqref{suppX}, there exists
$w_1\in \c_1$ such that $L_{w_1}$ is a subquotient of
$H^i\left(\fPhi\left(\bbK \otimes _{\O_\N}
\fPhi^{-1}(L_w)\right)\right)$ for some $i$. The object $\fPhi\left(
\bbK \otimes
_{\O_\N} \fPhi^{-1}(L_w) \right)$ is represented by the complex $$0\to
\Z(\Lambda^d(V)) *L_w \to \Z(\Lambda^{d-1}(V))*L_w \to \cdots
\to \Z(V)*L_w \to L_w \to 0.$$ But the Jordan-Hoelder series of
$Z_\lambda*L_w$ consists of $L_u$ with $u\LRleq w$. Hence
$\c_1\leq \c_2$. The Theorem is proved. \epf

\subsection{Proof of Proposition \ref{Duflo}}
The Proposition will be deduced from the next two Lemmas

\begin{Lem}\label{cha}
Let $j:O\imbed \N$ be an orbit of codimension $2m$,
 and  $\F\in D^\LG(\N)$
satisfy the following properties

i) $\F$ is an irreducible  perverse coherent sheaf with respect to the
perversity \eqref{pe}.

ii) $supp(\F)=\Obar$.

iii) $Hom_{D^{\LG}(\N)} (\O, \F[m])\ne 0$.

Then $\F\cong \Ohat_{\Obar}[-m]$.
\end{Lem}

\proof The condition
 $Hom_{D^{\LG}(\N)} (\O, \F[m])\ne 0$ is equivalent to the existence
of a nonzero $\LG$-invariant section of the coherent sheaf $H^{m}(\F)$
(where the cohomology is taken with respect to the
 usual $t$-structure on the derived category of coherent sheaves).
For a perverse coherent sheaf $\F$ on  $\Obar$ we have $H^i(\F)=0$ for $i<m$,
and
 $H^{m}(\F)$ is a torsion free sheaf on $\Obar$.
(Indeed, otherwise we would have a nonzero morphism defined on a
$\LG$-invariant open subscheme of $\Obar$ from
$V$ to $\F[i]$ where $i\leq m$ and $V$ is the nonderived
direct image of a vector
bundle under the locally 
closed imbedding of an orbit $O'\subset \Obar$, $O'\ne O$.
Since $V[-d]$ is a perverse coherent sheaf for $d=\frac{codim O'}{2}>m$
this would give an Ext of degree $i-d<0$ between perverse coherent sheaves,
which is impossible.)

 Thus  a nonzero
section of $H^m(\F)$ does not vanish on $O$. Also, $j^*(H^{m}(\F))$
 is an irreducible
$\LG$-equivariant vector bundle. Such a vector bundle
has a nonzero $\LG$-invariant section iff it is trivial; in which case
we have $\F\cong j_{!*}(\O_O[-m])\cong \Ohat_{\Obar}[-m],$
where the last equality is proved in  \cite{nilpokon}, Remark 11. \epf

\begin{Lem}\label{Homnenol}
We have $$Ext^{a(d_\c)}
_{\fPI}(L_0,L_{d_\c})
\ne 0,$$
where $a$ stands for Lusztig's $a$-function on $W$, see e.g.
\cite{cells2}, 1.1.
\end{Lem}

\proof The standard definition of a Duflo involution (see e.g.
\cite{cells2}, 1.3) shows that the costalk $j_e^!(L_d)$ has nonzero
cohomology in degree $a(d)$. We can think of $j_e^!(L_d)$
as an object in the $\bI$-equivariant derived category of $l$-adic
sheaves on the point. Moreover, it is a pull-back of an object in the
 $\bI$-equivariant derived category of $l$-adic
sheaves on the spectrum of a finite field.  The latter object is known to be
pure (cf. e.g. \cite{KGB}, Appendix, section A.7);
hence it is isomorphic to the
direct sum of its cohomology
(notice that Hom between two objects of
the bounded $\bI$-equivariant derived category of the point is identified
with Hom between corresponding complexes with constant cohomology
on $(\Pn)^{\rank(G)}$, $n\gg 0$.
Thus any pure object in the $\bI$-equivariant derived category of the point
 is isomorphic to the sum of its
cohomology by \cite{BBD}, Theorem 5.4.5).

 It follows that
$Hom _{\Db_I(\Fl)}(L_e,L_d[a(d)])\ne 0$. In view of Theorem 2 of \cite{1}
we will be done if we check that the map
$$Hom _{\Db_I(\Fl)}(L_e,L_d[a(d)])\to Hom_{\DIW}(\Delta_e, \Delta_e*L_d[a(d)]),
$$
sending $h$ to $id_{\Delta_e}*h$ is injective.

Recall (see e.g. \cite{Lco}) that $L_d[a(d)]$ is a direct summand in
$L_w*L_{w^{-1}}$ for any $w\in \c\cap \fW$ (e.g. for $w=d$).
Thus for any $h\in Hom_{D_I(\Fl)}(L_e,L_d[a(d)])$ the composition
\begin{equation}\label{co}
L_e\overset{h}{\To} L_d[a(d)]\to L_w*L_{w^{-1}}
\end{equation}
is nonzero for such $w$.

For $w\in W$ and $X,Y\in D(\Fl)$; or $X,Y\in D_I(\Fl)$
 we have a canonical isomorphism
\begin{equation}\label{i}
Hom(X*L_w,Y)\cong Hom(X,Y*L_{w^{-1}}).\end{equation}
In particular $Hom(L_e, L_w*L_{w^{-1}})\cong Hom(L_w,L_w)$
is a one dimensional space; thus multiplying
 $h\in Hom_{D_I(\Fl)}(L_e,L_d[a(d)])$ by a constant we can assume that
the composition \eqref{co} corresponds to $id\in Hom(L_w,L_w)$ under
the isomorphism \eqref{i}.
Then one can check that the composition
$$\Delta_e\overset{id_{\Delta_e}*h}{\To} \Delta_e* L_d[a(d)]
\to \Delta_e* L_w*L_{w^{-1}} $$
corresponds under \eqref{i} to $id\in Hom(\Delta_e*L_w, \Delta_e*L_w)$.
In particular, it is not equal to zero. \epf

\subsubsection{}%{Proof of Theorem \ref{Duflo} (conclusion)}
We are now ready to finish the proof of the Proposition. It
suffices to see that the object $(\fPhif)^{-1}(L_d)$ satisfies the
conditions of Lemma \ref{cha}. The first condition holds by
Theorem \ref{tstr}. The second one holds by Theorem \ref{kle}(a).
 Finally, to check
condition (iii) notice that by Lemma \ref{Homnenol} we have
$$\begin{array}{ll}
Hom(\O_\N,\fPhif^{-1}(L_d)[a(d_\c)])=Hom(\O_\N, p_{Spr*}(\fPhi^{-1}
(L_d)[a(d_\c)])\\ = Hom(\O_\Nt, \fPhi^{-1}(L_d)[a(d_\c)])=
Hom(L_e,L_d[a(d_\c)])
\ne 0.
\end{array}
$$
By  \cite{cells4}, Theorem 4.8(c) we have
$a(d_\c)= \frac{\codim (\LG(N_\c))}{2}$, which implies condition (iii). \epf

\subsection{Proof of Proposition \ref{risPhi}}\label{proofrisPhi}
 If $\rho$ is trivial then
\eqref{onorbit} follows from Proposition \ref{Duflo}. Applying
\eqref{rc} we see that  \eqref{onorbit} holds when
$\rho=Res^\LG_{Z_\LG(N_\c)}(V)$ for $V\in Rep(\LG)$.

Let now $\rho$ be arbitrary. Let $\L\in Coh^\LG(\N)$ be some sheaf
supported on the closure of $\LG(N_\c)$, and such that
$\L|_{\LG(N_\c)}\cong \L_\rho$.
 We can choose a short exact sequence
\begin{equation}\label{eses} W\otimes \O \overset{\phi}{\To}
V\otimes \O\to \L\to 0, \end{equation}
 $V, W\in Rep(\LG)$.

Then we get an exact sequence
$$ W|_{Z_\LG(N_\c)}\to V|_{Z_\LG(N_\c)} \to \rho \to 0$$
in $Rep(Z_\c)$, and hence an exact sequence in $\fPI^\c$:
$$r_\c(W|_{Z_\c})\circ L_{d_\c}\to
r_c(V|_{Z_\c})\circ L_{d_\c} \to r_c(\rho)\circ L_{d_\c}
\to 0;$$ by \eqref{rc} it can be written as
\begin{equation}\label{ese}
\Z(W)* L_{d_\c}\to \Z(V)* L_{d_\c} \to r_c(\rho)\circ L_{d_\c} \to
0.
\end{equation}

On the other hand, consider the tensor product of \eqref{eses} by
$j_*(\O)$ where $j$ stands for the imbedding $\LG(N_\c)\imbed \N$
(and $j_*$ is the non-derived direct image). We get a short exact
sequence
\begin{equation}\label{esess}
W\otimes j_*(\O) \to V\otimes j_*(\O) \to L' \to 0
\end{equation}
where $L'|_{\LG(N_\c)}\cong \L_\rho$. Theorem \ref{tstr} and the
definition of a perverse coherent sheaf show that the functor
$\F\mapsto \Phi_\c(\F)[-m]$ is exact with respect to the {\it
standard} $t$-structure on the category $D^\LG_{N_\c}(\N)$.
Applying this functor to \eqref{esess} we get an exact sequence in
$\fPfI^\c$, which by Theorem \ref{Duflo} has the form
\begin{equation}\label{seses}
\Z(W)* L_{d_\c} \to \Z(V)* L_{d_\c} \to \Phi_\c(\L_\rho[-m]) \to
0.
\end{equation}

Lemma \ref{compaN} implies that \eqref{seses} is isomorphic to
\eqref{ese} (or rather to its image in the quotient category
$\fPfI^\c$); in particular, \eqref{onorbit} holds. \epf

\subsubsection{Proof of Corollary \ref{HZ}} 
The functor $\Phi_\c$ is an equivalence, thus Proposition 
 \ref{risPhi} implies that the functor $r_\c^f: Rep(Z_\c)\to A_\c^f\cong 
Rep(H_\c)$
is fully faithful. The functor of restriction of a representation to 
a subgroup can only be fully faithful if the subgroup coincides with
the whole group, thus $H_\c=Z_\c$. \epf

\footnotesize{
Department of Mathematics,
 Massachusetts Institute of Technology,
Cambridge MA,
02139, USA;\\ 
%\hphantom{x}\ab\, 
{\tt bezrukav@math.mit.edu}}

\end{document}